\documentclass[11pt,a4paper]{article}
\usepackage[T1]{fontenc}
\usepackage{lmodern,amsmath,amssymb,amsthm,mathtools,booktabs,array}
\usepackage[margin=26mm]{geometry}
\usepackage{microtype}
\usepackage[hidelinks]{hyperref}
\usepackage{enumitem}
\setlist{itemsep=3pt,topsep=5pt}
\newtheorem{theorem}{Theorem}[section]
\newtheorem{lemma}[theorem]{Lemma}
\newtheorem{proposition}[theorem]{Proposition}
\newtheorem{corollary}[theorem]{Corollary}
\theoremstyle{remark}

\newcommand{\Fp}{\mathbb F_p}
\newcommand{\ind}{\mathbf1}
\title{The strong cylinder conjecture over prime fields}
\date{}
\author{Peter Sziklai}
\begin{document}
\maketitle
\begin{abstract}
Let $p$ be prime and let $S\subseteq\mathbb F_p^3$ have $p^2$ points,
with every affine plane meeting $S$ in a multiple of $p$ points.
We give a self-contained argument that $S$ is a union of $p$ parallel
lines.
\end{abstract}

\section{Statement and proof structure}

\begin{theorem}[Strong cylinder conjecture over prime fields]
\label{thm:prime-strong-cylinder}
Let $p$ be prime. If $S\subseteq\mathbb F_p^3$ has $p^2$ points
and every affine plane meets $S$ in a multiple of $p$ points, then
$S$ is a union of $p$ distinct parallel affine lines.
\end{theorem}

Call a set satisfying the hypotheses \emph{admissible}, and call the
asserted union a \emph{cylinder}. A plane containing more than $p$
points of $S$ is \emph{heavy}; a plane containing exactly $p$ points
is \emph{light}. A line containing at least two points of $S$ is a
\emph{secant}, and an $m$-secant contains exactly $m$ points.
The prime-field conjecture and its relation to divisible codes and
signed cylinder decompositions are discussed in \cite{KM,KMNS}.
The conjecture is due to Ball \cite{Ball}, and it is known to hold
for $p\le7$ \cite{DDMS,KM}.
The prime-field assumption is essential; this manuscript does not
assert the analogous statement over all finite fields.

The elementary part first reduces the problem to constructing, for a
hypothetical noncylinder, canonical linear forms $\lambda_P$ such that
\begin{equation}\label{eq:roadmap-harmonic}
 \lambda_{P_i}(P_i)=-2,\qquad
 \lambda_{P_i}(v)=\frac{2}{m-1}
          \sum_{j\ne i}\frac1{t_j-t_i}
 \quad(P_i=P_0+t_i v\in S\cap\ell).
\end{equation}
Here points are represented by $(1,P)\in\Fp^4$, and the direction
vector $v$ has first coordinate zero. On a light plane $H$, the map
\[
 A_H=\sum_{Q\in S\cap H}Q\,\lambda_Q|_{\widetilde H}
\]
acts on its three-dimensional vector space. Formula
\eqref{eq:roadmap-harmonic} makes each point of $S\cap H$ an eigenvector.
Elementary geometry forces $A_H$ to be scalar, and its trace forces it
to be zero. This makes the light section a full line.

Sections~\ref{sec:chow-symbolic}--\ref{sec:harmonic} construct the forms
and prove \eqref{eq:roadmap-harmonic}; Section~\ref{sec:completion}
gives the final argument. Throughout the algebraic construction
$p\ge5$ and $N=p-1$.

\paragraph{Coordinate conventions.}
Physical projective points are columns $P=(1,p_1,p_2,p_3)$ in a vector
space $V=\Fp^4$; its affine chart is fixed by $w(P)=P_0=1$.
The variables $X,T$ are dual coordinates, so a physical point gives
the linear form $L_P(X)=P^{\mathsf T}X$. A decomposable dual bivector
$X\wedge T$ represents the physical line annihilated by $X$ and $T$.
The auxiliary form $G(z)$ below is in physical coordinates $z$.
Consequently, physical lines contained in a plane correspond on the
Klein quadric to dual two-spaces containing that plane's normal.
All physical point-centered expansions used to define $\lambda_P$
preserve the fixed affine normalization $w=1$. The dual normalizations
in the heavy-plane calculation may use arbitrary rational projective
coordinates; the truncated construction is equivariant under these
changes, as proved below.

\section{Elementary geometry}

\begin{lemma}[A full line forces a cylinder]\label{lem:full-line}
If an admissible set contains a full affine line, it is a cylinder.
\end{lemma}
\begin{proof}
Let the contained line have direction $v$. Every plane transverse to $v$
meets this line, and therefore contains at least $p$ points of $S$.
The $p$ planes in each such parallel class partition $S$, so each contains
exactly $p$ points.
Fix $P\in\Fp^3$, and put $r=|S\cap(P+\Fp v)|$.
There are $p^2$ planes through $P$ transverse to $v$. A point $Q$ outside
$P+\Fp v$ lies in $p$ of them; a different point of that line lies in none.
Counting incidences gives
\[
 p^3=p^2\ind_S(P)+p(p^2-r),\qquad r=p\ind_S(P).
\]
Thus every line parallel to $v$ is empty or full, proving the assertion.
\end{proof}

\begin{lemma}[Local positive weights]\label{lem:local-weights}
For $P\in S$ and planes $H\ni P$, put $c_P(H)=|S\cap H|/p-1$.
Then these are nonnegative integers and
\begin{equation}\label{eq:local-weights}
 \sum_{H\ni P}c_P(H)=p-1,
 \qquad
 \sum_{H\supseteq\ell}c_P(H)=|S\cap\ell|-1
 \quad(\ell\ni P).
\end{equation}
In particular, $S$ is contained in the union of the heavy planes through $P$.
\end{lemma}
\begin{proof}
There are $p^2+p+1$ planes through $P$, and each other point lies in $p+1$
of them. This gives the first identity. The $p+1$ planes through a line
$\ell$ have total intersection size $p^2+p|S\cap\ell|$, giving the second.
If $Q\in S\setminus\{P\}$, the latter sum for $PQ$ is positive.
\end{proof}

\begin{lemma}[$(p-1)$ points on a line suffice]
\label{lem:almost-full-line}
Let $p\ge3$.  If an admissible set has exactly $p-1$ points on an affine
line, it contains a full affine line and is a cylinder.
\end{lemma}
\begin{proof}
Let $\ell$ contain $p-1$ points of $S$.  By
Lemma~\ref{lem:local-weights}, the heavy weights in the pencil
of planes containing $\ell$ sum to $p-2$.  Of its $p+1$ planes,
at least three consequently meet $S$ in exactly $p$ points.
For each such plane $\pi$ write
\[
 S\cap\pi=(S\cap\ell)\cup\{Q_\pi\},\qquad Q_\pi\notin\ell.
\]

For $P\in S\cap\ell$, the total heavy weight through $P$ is
$p-1$, of which $p-2$ is contributed by planes containing
$\ell$.  Hence there is a unique further heavy plane $H_P$;
it has weight one and does not contain $\ell$.
The line $PQ_\pi$ is a two-secant.  Its heavy-pencil weight is
therefore one.  A plane containing both $PQ_\pi$ and $\ell$
would be the light plane $\pi$, so the required heavy plane is
$H_P$.  Thus every $Q_\pi$ belongs to every $H_P$.

Choose two distinct light planes in this pencil.  Their points
$Q_1,Q_2$ are distinct, and the line $m=Q_1Q_2$ lies in every
$H_P$.  The $p-1$ planes $H_P$ are distinct: a plane containing
two different points of $S\cap\ell$ would contain $\ell$.
Applying the radial identity at $Q_1$ now gives
\[
 |S\cap m|-1
 =\sum_{H\supseteq m}\left(\frac{|S\cap H|}{p}-1\right)
 \ge p-1.
\]
Since $m$ is an affine line, it is full.  Lemma~\ref{lem:full-line}
finishes the proof.
\end{proof}

\section{The R\'edei product and Frobenius bilinears}\label{sec:chow-symbolic}

Set
\[
 B_{ij}=X_iX_j^p-X_jX_i^p\qquad(0\le i<j\le3).
\]
The multiplicity of a homogeneous polynomial at a projective point means
membership in the corresponding power of the homogeneous point ideal.

\begin{lemma}\label{lem:translation-invariants-three-pairs}
Let $x=(x_1,x_2,x_3)$ and $t=(t_1,t_2,t_3)$.  A polynomial
$F\in \Fp[x,t]$ satisfying
\[
 F(x+\lambda t,t)=F(x,t)\quad\text{in }\Fp[x,t,\lambda]
\]
belongs to the subalgebra generated by the $t_i$ and the three determinants
$u_{ij}=x_it_j-x_jt_i$.
\end{lemma}

\begin{proof}
Consider
\[
 \mathcal R=\Fp[t_1,t_2,t_3,u_{12},u_{13},u_{23}]/
       (t_1u_{23}-t_2u_{13}+t_3u_{12}).
\]
The displayed relation is irreducible, and the indicated substitution
embeds $\mathcal R$ into $\Fp[x,t]$: this is immediate after inverting $t_1$, where
$u_{23}=(t_2u_{13}-t_3u_{12})/t_1$ and
$t_1,t_2,t_3,u_{12},u_{13}$ are algebraically independent.
After the same localization,
\[
 \Fp[x,t,t_1^{-1}]
 =\Fp[t_1^{\pm1},t_2,t_3,u_{12},u_{13},x_1],
\]
and translation changes only $x_1$, by $x_1\mapsto x_1+\lambda t_1$.
A polynomial invariant under this substitution is independent of $x_1$.
Consequently $F\in \mathcal R[t_1^{-1}]$.

It remains to cancel possible powers of $t_1$ in the denominator.
The induced map
\[
 \mathcal R/(t_1)\longrightarrow \Fp[x,t]/(t_1)
\]
is injective.  Indeed, its source is the domain
\[
 \Fp[t_2,t_3,u_{12},u_{13},u_{23}]/(t_2u_{13}-t_3u_{12}),
\]
and after inverting $t_2$ injectivity follows from the algebraic
independence of $t_2,t_3,u_{12}=t_2x_1$ and
$u_{23}=t_3x_2-t_2x_3$.  Thus any element of $\mathcal R$ which is divisible by
$t_1$ in the ambient polynomial ring is already in $t_1\mathcal R$.
Writing $F=P/t_1^N$ and applying this observation repeatedly proves
$F\in\mathcal R$.
\end{proof}

\begin{theorem}\label{thm:minimal-symbolic-piece}
Let $1\le r\le p-1$. Suppose that $H$ is a homogeneous polynomial
in four variables over $\Fp$, of degree $r(p+1)$ and multiplicity at least
$r$ at every point of $\mathrm{PG}(3,p)$.  Then $H$ is an $\Fp$-linear
combination of products of $r$ of the polynomials $B_{ij}$.

There is no nonzero homogeneous polynomial of degree less than
$r(p+1)$ with these multiplicities.
\end{theorem}

\begin{proof}
We first prove the assertion at degree $r(p+1)$.  On the chart $X_0=1$
write $x_i=X_i/X_0$, $h(x)=H(1,x)$ and $U_i=x_i^p-x_i$.
The ideals of distinct affine rational points are pairwise comaximal.
Their intersection is $J=(U_1,U_2,U_3)$, so the multiplicity hypothesis
gives $h\in J^r$.

The polynomial ring $\Fp[x_1,x_2,x_3]$ is free over
$\Fp[U_1,U_2,U_3]$ with basis $x^\beta$, $0\le\beta_i<p$.
Moreover, the leading monomials of $U^\alpha x^\beta$ are all distinct,
and their degrees are $p|\alpha|+|\beta|$.  Since
\[
 p(r+1)>r(p+1),
\]
the unique expansion of $h\in J^r$ contains only terms with
$|\alpha|=r$.  Thus
\[
 h(x)=\sum_{|\alpha|=r} A_\alpha(x)U^\alpha,
 \qquad \deg A_\alpha\le r.
\]
Homogenize each coefficient to degree $r$, and set
\[
 V_i=X_i^p-X_iX_0^{p-1},\qquad
 A(X,T)=\sum_{|\alpha|=r}A_\alpha^{\mathrm h}(X)T^\alpha.
\]
Then $A$ has bidegree $(r,r)$ in $X=(X_0,X_1,X_2,X_3)$ and
$T=(T_1,T_2,T_3)$, respectively, and
\[
 H(X)=A(X,V(X)).
\]

Fix a rational point at infinity $P=(0,t_1,t_2,t_3)$.
In local increments at $P$, each $V_i$ is congruent to the constant
$t_i$ modulo terms of order at least $p-1$, hence modulo terms of
order at least $r$.  For every multi-index $\beta$ with
$|\beta|=r-1$, the multiplicity hypothesis therefore gives
\[
 K_\beta(t):=(\partial_X^\beta A)(0,t;t)=0.
\]
The polynomial $K_\beta(T)$ is homogeneous of degree $r+1\le p$.
No nonzero homogeneous polynomial of degree at most $p$ vanishes at
every point of $\mathrm{PG}(2,p)$: its restriction to each rational
line would have $p+1$ zeros and hence vanish identically, forcing more
linear factors than its degree.  Consequently every $K_\beta$ is
identically zero.

Put $\mathcal D_T=\sum_{i=1}^3T_i\partial_{X_i}$.
Since $\partial_X^\beta A$ is linear in $X$ when $|\beta|=r-1$,
\[
 \partial_X^\beta(\mathcal D_TA)
 =\mathcal D_T(\partial_X^\beta A)
 =(\partial_X^\beta A)(0,T;T)=K_\beta(T)=0.
\]
The polynomial $\mathcal D_TA$ is homogeneous of $X$-degree $r-1<p$;
all derivatives of that order vanish, so $\mathcal D_TA=0$.
Taylor expansion is valid with invertible factorials through degree
$r<p$, and yields the polynomial identity
\[
 A(X_0,X_1+\lambda T_1,X_2+\lambda T_2,X_3+\lambda T_3;T)
 =A(X;T).
\]

Write $A=\sum_{s=0}^rX_0^sA_s(X_1,X_2,X_3;T)$.
Each $A_s$ has bidegree $(r-s,r)$ and is translation invariant.
By Lemma~\ref{lem:translation-invariants-three-pairs}, it is a
polynomial in the $T_i$ and $X_iT_j-X_jT_i$.
The invariant-ring presentation is bihomogeneous; taking its
$(r-s,r)$ component shows that every term uses exactly $r-s$
determinants and $s$ additional factors $T_i$.
Thus $A$ is a linear combination of products of $r$ factors chosen
from
\[
 X_0T_i\quad(1\le i\le3),\qquad
 X_iT_j-X_jT_i\quad(1\le i<j\le3).
\]
Substitute $T_i=V_i(X)$.  These factors become, respectively,
$B_{0i}$ and $B_{ij}$, proving the assertion.

For the final assertion, suppose that $0\ne H$ has degree
$d<r(p+1)$ and the stated multiplicities.  If $d<pr$, the same free
module expansion on $X_0=1$ contradicts $h\in J^r$ immediately.
Otherwise put $s=d-pr$, so $0\le s<r$.  The expansion again has only
$|\alpha|=r$, now with coefficient degree at most $s$.
After homogenizing those coefficients to degree $s$, write
$H=A(X,V)$ with bidegree $(s,r)$.
For every $|\beta|=s$, multiplicity at infinity implies
$(\partial_X^\beta A)(0,t;t)=0$ for all rational $t$.
This is a homogeneous polynomial of degree $r\le p-1$ in $t$, hence
identically zero.  Since these derivatives are independent of $X$
and $s!$ is invertible, every coefficient of $A$ vanishes, a
contradiction.
\end{proof}

The \emph{Hasse Taylor term of order $r$} of a polynomial $F$ at $A$
is the homogeneous degree-$r$ part in $Z$ of $F(A+Z)$, defined
without division by factorials.

\begin{corollary}\label{cor:chow-partials-bilinears}
Let $R$ be a homogeneous polynomial of degree $p^2$ in four variables.
Suppose that at every rational projective point all its Hasse Taylor
terms of orders $1,\ldots,p-1$ vanish.  Then every ordinary partial
derivative of $R$ is a homogeneous polynomial of degree $p-1$ in the
six Frobenius bilinears $B_{ij}$.
\end{corollary}

\begin{proof}
Each partial derivative has degree $p^2-1=(p-1)(p+1)$ and
multiplicity at least $p-1$ at every rational projective point.
Apply Theorem~\ref{thm:minimal-symbolic-piece} with $r=p-1$.
\end{proof}

\begin{lemma}\label{lem:admissible-chow-jets}
Suppose that $S\subseteq\mathbb F_p^3$ has $p^2$ points and every
affine plane meets $S$ in a multiple of $p$ points.  Put
\[
 L_P(X)=X_0+P_1X_1+P_2X_2+P_3X_3,
 \qquad R(X)=\prod_{P\in S}L_P(X).
\]
We call $R$ the \emph{R\'edei product} of $S$.
At every rational point $A$ of $\mathrm{PG}(3,p)$, all Hasse Taylor
terms of $R(A+Z)$ of orders $1,\ldots,p-1$ vanish.
Consequently every partial derivative of $R$ is a homogeneous
polynomial of degree $p-1$ in the six $B_{ij}$.
\end{lemma}

\begin{proof}
For every $A\in\mathbb F_p^4$ the plane-divisibility assumption,
including the case of a constant linear expression, gives
\[
 \sum_{P\in S}L_P(A)^{p-1}=0.
\]
The left side is a polynomial of total degree $p-1$, hence degree
less than $p$ in each variable, and therefore is identically zero.
All multinomial coefficients in this degree are nonzero in
$\mathbb F_p$.  Comparing coefficients shows that the sum over
$(1,P)$, $P\in S$, of every homogeneous polynomial of degree $p-1$
is zero.

If some $L_P(A)$ vanishes, the corresponding affine plane contains
at least $p$ points of $S$.  The product $R(A+Z)$ then has at least
$p$ linear factors vanishing at $Z=0$, so it has no terms of degree
less than $p$.

Otherwise all $L_P(A)$ are nonzero.  For $1\le j\le p-1$,
\[
 \sum_{P\in S}\left(\frac{L_P(Z)}{L_P(A)}\right)^j
 =\sum_{P\in S}L_P(Z)^jL_P(A)^{p-1-j}=0.
\]
The last equality is the homogeneous moment identity just proved.
Newton's identities, with $1,\ldots,p-1$ invertible, imply that the
elementary symmetric polynomials of these ratios vanish in degrees
$1,\ldots,p-1$.  Thus
\[
 R(A+Z)=R(A)+\text{terms of total $Z$-degree at least $p$}.
\]
The final conclusion follows from
Corollary~\ref{cor:chow-partials-bilinears}.
\end{proof}

\section{The R\'edei tensor and its complementary form}\label{sec:fundamental-lines}

The following auxiliary polynomial transfers the identities for the partial
derivatives of $R$ to the original projective space. Put $N=p-1$ and
assume $p$ is odd. Write $R_i=\partial R/\partial X_i$.
Recall the proved representation
\[
 R_i(X)=\Phi_i(X\wedge X^{[p]}),\qquad \deg\Phi_i=N.
\]
For independent four-vectors $X,T$, write
$\Psi_i(X,T)=\Phi_i(X\wedge T)$. Each $\Psi_i$ has bidegree $(N,N)$.
The polynomial $\Phi_i$ is determined only modulo the Pl\"ucker relation
$B_{01}B_{23}-B_{02}B_{13}+B_{03}B_{12}$, but $\Psi_i$ is well defined.
Indeed, the Pl\"ucker relation vanishes on $X\wedge T$, and if two
choices of $\Phi_i$ agree after the substitution $T=X^{[p]}$, the
difference of the corresponding $\Psi_i$ is a form of bidegree $(N,N)$
killed by this substitution. The substitution is injective on such
forms: every $X$-exponent is at most $N<p$, so the exponents
$\alpha+p\beta$ determine $\alpha$ and $\beta$.
Since $N=p-1$ is even, $\Psi_i(T,X)=\Psi_i(X,T)$.

\begin{lemma}\label{lem:fundamental-form}
There is a nonzero reduced homogeneous polynomial
$G(z_0,z_1,z_2,z_3)$ of degree $2p-3$, with each variable degree at
most $p-1$, such that
\[
 S\subseteq\mathcal T(G):=
 \{[P]\in\operatorname{PG}(3,p):D_P^{p-1}G=0\},
 \qquad D_P=\sum_iP_i\partial_{z_i}.
\]
Here the equation $D_P^{p-1}G=0$ is a polynomial identity, not only
vanishing at one point. The construction is compatible with projective
$\mathbb F_p$-coordinate changes in the truncated algebra
$\mathcal A=\mathbb F_p[z_0,z_1,z_2,z_3]/(z_0^p,z_1^p,z_2^p,z_3^p)$.
\end{lemma}

\begin{proof}
First lift Euler's identity to
\[
 \sum_iX_i\Psi_i(X,T)=\sum_iT_i\Psi_i(X,T)=0.
\]
For the first identity, the left side has bidegree $(p,p-1)$ and
vanishes after $T=X^{[p]}$. Monomials can collide under this
substitution only when their $X$-parts are pure $p$th powers.
Consequently the left side has the form $\sum_iX_i^pC_i(T)$,
with $\deg C_i=p-1$. Substituting $T=X$ makes it zero because
$X\wedge X=0$. The monomial supports of $X_i^pC_i(X)$ are disjoint,
so every $C_i=0$. Interchanging $X$ and $T$ gives the other identity.

Equality of mixed partial derivatives of $R$ similarly lifts to
\[
 \partial_{X_j}\Psi_i=\partial_{X_i}\Psi_j,
 \qquad
 \partial_{T_j}\Psi_i=\partial_{T_i}\Psi_j.
\]
Indeed, substitution $T=X^{[p]}$ is injective on the relevant
polynomials of $X$-degree $p-2<p$; the second system follows by
interchanging the two vectors.

Polarize the two blocks of size $N$, using that $N!\ne0$. There is
a tensor $C$ with
\[
 \Psi_i(X,T)=C(e_i,X^N,T^N).
\]
The two mixed-partial systems allow the output index to exchange with
any index in either block, so $C$ is fully symmetric of order
$2N+1=2p-1$. We call $C$ the \emph{R\'edei tensor}: it is a
symmetric multilinear form in $2p-1$ arguments, or equivalently
its coefficient array, unchanged by permuting the indices.
Here $X^N$ and $T^N$ denote $N$ repeated arguments $X$ and $T$.
The horizontal identity is $C(X^p,T^N)=0$.
In characteristic $p$, mixed multinomial coefficients in the block
of length $p$ vanish. Polarizing the remaining block shows that
\[
 C_\gamma=0\quad\hbox{if any }\gamma_i\ge p,
 \qquad |\gamma|=2p-1.
\]
Thus $C$ is a linear functional on $\mathcal A_{2p-1}$.
The perfect multiplication pairing into the one-dimensional socle
$\mathcal A_{4p-4}$ identifies it with the nonzero form
\[
 G=\sum_{\substack{|\gamma|=2p-1\\0\le\gamma_i\le N}}C_\gamma
     z^{N\mathbf1-\gamma}\in\mathcal A_{2p-3}.
\]
Nonzeroness follows from $\nabla R\ne0$, since $R$ is squarefree.
The pairing is equivariant up to the scalar action on the socle.
Explicitly, let the chosen vector representatives transform as $P'=MP$,
without further individual rescaling. Then
\[
 G'(z)=\det(M)^N\operatorname{trunc}\bigl(G(M^{-1}z)\bigr).
\]
Since $\det(M)\in\Fp^\times$ and $N=p-1$, we have $\det(M)^N=1$,
so $G$ in fact transforms exactly, not merely up to a scalar.
Here $\operatorname{trunc}$ discards monomials containing $z_i^p$;
it is reduction modulo the ideal of pure $p$th powers. Adding a term
$z_i^pH_i$ of degree $2p-3$ cannot change any derivative $D_P^NG$,
since $\deg H_i=p-3<N$. Consequently
\[
 D_{MP}^NG'(z)=\det(M)^N(D_P^NG)(M^{-1}z),
\]
as an ordinary polynomial identity. The output degree is $p-2<p$.

To check the asserted condition at $P\in S$, choose rational
coordinates with $P=e_0$. On $X_0=0$, the squarefree product gives
$R_i=0$ for $i\ne0$. Substitution injectivity, now with $X$-degree
$N<p$, gives $\Psi_i(X,T)=0$ for $i\ne0$ whenever $X_0=T_0=0$.
Polarization says precisely that $C_\gamma=0$ when $\gamma_0=0$.
Hence $G$ has no monomial with $z_0$-exponent $N$, or equivalently
$\partial_{z_0}^NG=0$. Equivariance gives $D_P^NG=0$ in the
original coordinates. The ideal of $p$th powers is stable under
linear coordinate changes and ordinary derivatives, which justifies
working throughout with reduced representatives.
\end{proof}

\section{Local products and canonical endpoint planes}

We retain $N=p-1$, the universal gradient
$\nabla R(X)=\Phi(X\wedge X^{[p]})$, and its complementary form $G$.
We first supply the normal-order and local-product arguments used
below.

\begin{lemma}\label{lem:heavy-klein-normal-order}
If a physical plane $H$ contains $kp$ points of $S$, where $k\ge1$,
the universal section $\Phi$ vanishes to normal order exactly $k-1$
along the Klein plane $\beta_H$ of physical lines contained in $H$.
Here the normal order of $\Phi=(\Phi_0,\ldots,\Phi_3)$ is the minimum
of the orders of its components: every $\Phi_i$ vanishes to order at
least $k-1$ along $\beta_H$, with equality for at least one $i$.
\end{lemma}
\begin{proof}
Put the dual point of $H$ at $e_0$ and write
\[
 a=(B_{01},B_{02},B_{03}),\qquad
 b=(B_{23},-B_{13},B_{12}).
\]
The Klein coordinate ring is $\mathbb F_p[a,b]/(a\cdot b)$, and $\beta_H$
has ideal $(b)$. Its relation is bihomogeneous. Consequently a
homogeneous component $\Phi_i$ of degree $N$ has a unique direct-sum
decomposition by $b$-degree, with the term of $b$-degree $j$ a class
$A_j$ of bidegree $(N-j,j)$. The least nonzero $j$ is precisely its
order along $\beta_H$.

The substitution $b=a\times T$ is injective on the Klein ring:
over $a\ne0$, every $b$ with $a\cdot b=0$ has this form. If $A_j$
is nonzero, $A_j(a,a\times T)$ is therefore nonzero and has bidegree
$(N,j)$ in $(a,T)$. Substituting $T=a^{[p]}$ remains injective on
these forms. Indeed, every exponent of $a$ in its first block is at
most $N<p$, so the resulting exponents $\alpha+p\beta$ determine
both multiindices uniquely.

On the dual affine chart $X=(1,x)$, use the formal parameters
$u_i=x_i^p-x_i$. They are regular local parameters at $x=0$, with
\[
 x_i=-u_i-u_i^p-u_i^{p^2}-\cdots.
\]
Thus $a=u$ and
\[
 b=u\times u^{[p]}+
        \text{terms of degree at least }p^2+1.
\]
The first nonzero term contributed by $A_j$ has ordinary degree
$(N-j)+(p+1)j=N+pj$, and is nonzero by the preceding injectivity.
Different $j$ give different first degrees. Hence the local order
of $\Phi_i(X\wedge X^{[p]})$ is $N+pj$.

Exactly $kp$ factors of the squarefree R\'edei product vanish at the
dual point. Their distinct linear initial forms have squarefree
product of degree $kp$, whose gradient is not identically zero.
The minimum local order among the partials of $R$ is consequently
$kp-1=N+p(k-1)$. Comparing the two descriptions proves the claim.
\end{proof}

\begin{lemma}\label{lem:local-klein-heavy-product}
Let $P=e_0\in S$. On the Klein plane $\alpha_P$ of physical lines
through $P$, use $b=(B_{23},-B_{13},B_{12})$. Then
\[
 \Phi|_{\alpha_P}=e_0h(b),\qquad
 h(b)=c_P\prod_{H\ni P,\ H\text{ heavy}}
                 \ell_H(b)^{|S\cap H|/p-1},\quad c_P\ne0,
\]
where $\ell_H$ is the equation of $H$ in the three direction
coordinates at $P$. With the same tensor normalization and oriented
coordinates, $[z_0^{p-2}]G=h(z_1,z_2,z_3)$.
\end{lemma}
\begin{proof}
The fundamental-point condition gives the first displayed form.
On $X_0=0$, the partial derivative $R_0=\partial R/\partial X_0$ is the product of the restrictions
of all factors other than $L_P$. Grouping those factors by physical
lines through $P$ gives, up to a nonzero scalar,
\[
 \left.R_0\right|_{X_0=0}
   =\prod_{d\in\operatorname{PG}(2,p)}(d\cdot x)^{m_d-1},
 \qquad m_d=|S\cap(P+\mathbb F_p d)|.
\]
Here the affine-line notation means the original affine line through
$P$ of projective direction $d$. The radial identity is
$m_d-1=\sum_{H\supset Pd}(|S\cap H|/p-1)$. Interchanging the products
therefore expresses the last polynomial as
\[
 c_P\prod_{H\ni P,\ H\text{ heavy}}
       \left(\prod_{d\in\operatorname{PG}(H/P)}d\cdot x\right)
                  ^{|S\cap H|/p-1}.
\]
For each rational two-dimensional subspace of directions, the product
of its $p+1$ linear forms is a nonzero scalar multiple of the relevant
linear form in $x\times x^{[p]}$. This is the usual binary identity
$\prod_{[a:b]\in\operatorname{PG}(1,p)}(ax+by)
 \doteq x^py-xy^p$ after a rational change of basis.
Thus the claimed product for $h$ holds after substituting
$b=x\times x^{[p]}$. This substitution is injective on degree-$N$
forms: first substitute $b=x\times T$, which is dominant, and then
$T=x^{[p]}$, which is injective because the first-block degree is
$N<p$. This proves the product identity itself.

For completeness, write
$h(b)=\sum_{a+b+c=N}h_{abc}b_1^ab_2^bb_3^c$ and specialize the two
ternary blocks to $X=(x_1,x_2,0)$ and $T=(t_1,0,t_3)$.
Their oriented wedge is $(x_2t_3,-x_1t_3,-x_2t_1)$. The indicated
monomial contributes
\[
 (-1)^{b+c}h_{abc}
   x_1^b x_2^{a+c}t_1^c t_3^{a+b}.
\]
In $C(e_0,X^N,T^N)$ the polarization coefficient of this monomial
is $\binom Nb\binom Nc=(-1)^{b+c}$ times
$C_{(1,b+c,a+c,a+b)}$. Thus this tensor entry equals $h_{abc}$.
Complementing its four exponents gives
$(N-1,a,b,c)$, exactly the coefficient of
$z_0^{p-2}z_1^az_2^bz_3^c$ in $G$. This proves the final assertion.
\end{proof}

\begin{lemma}\label{lem:heavy-truncated-power}
Let $H$ be a heavy plane with $|S\cap H|=(r+1)p$. Then its linear equation $\ell_H$ satisfies
\[
 G\in(\ell_H^r)
 \quad\hbox{in }\mathbb F_p[z_0,z_1,z_2,z_3]/(z_0^p,z_1^p,z_2^p,z_3^p).
\]
\end{lemma}
\begin{proof}
Choose coordinates with $\ell_H=z_3$. In the universal expression
$\Phi_i(X\wedge T)=C(e_i,X^N,T^N)$, put
$X=e_3+U$, with $U_3=0$, and $T_3=0$. The three Klein coordinates
normal to the plane of lines in $H$ are linear in $U$. By
Lemma~\ref{lem:heavy-klein-normal-order}, every component $\Phi_i$
vanishes to order at least $r$ along this plane; only this lower
bound is used here. It therefore kills every coefficient of
$U$-degree $j<r$.
If a nonzero tensor index $\gamma$ has
$\gamma_3=N-j$ with $0\le j<r$, choose $i\ne3$ with $\gamma_i>0$
and split the other indices after removing $e_i$ into multiindices
$\alpha,\beta$ of sizes $j,N$. The coefficient of $U^\alpha T^\beta$
is a nonzero factorial multiple of $C_\gamma$; all relevant
factorials are below $p$. Thus $C_\gamma=0$ whenever
$\gamma_3\ge N-r+1$. The complementary exponent in $G$ is
$N-\gamma_3$, so the reduced representative is divisible by $z_3^r$.
This is exactly the stated divisibility in the truncated algebra.
\end{proof}

\begin{lemma}\label{lem:degree-p-gluing}
Let $\ell_1,\ldots,\ell_h$ be pairwise non-proportional ternary linear forms whose
span has dimension three. Suppose positive integers $r_i$ sum to
$p-1$, and put $H=c\prod_i\ell_i^{r_i}$ with $c\ne0$.
If a homogeneous degree-$p$ form $K$ satisfies
\[
 K\in(\ell_i^{r_i})+(z_1^p,z_2^p,z_3^p)
 \quad\hbox{for every }i,
\]
then $K=B^p+LH$ for linear forms $B,L$ over $\Fp$.
\end{lemma}
\begin{proof}
Write $K-B_i^p\in(\ell_i^{r_i})$ for suitable linear forms $B_i$.
Choose three independent $\ell_i$ as coordinate forms. The
restrictions of $K$ to their three coordinate lines are $p$th powers
of linear forms, and agreement at the three vertices gives a single
linear $B$ whose $p$th power agrees with all three restrictions.
For any other $\ell_i$, the restriction of $K-B^p$ to its projective
line is a $p$th power of a linear form. It vanishes at at least two
distinct intersections with the three coordinate lines, and is
therefore identically zero. Consequently $B_i-B$ is divisible by
$\ell_i$. Since $r_i<p$, the difference $(B_i-B)^p$ is divisible by
$\ell_i^{r_i}$, and hence so is $K-B^p$. The distinct factors are
coprime, so $H$ divides $K-B^p$. The remaining quotient has degree one.
\end{proof}

\begin{proposition}\label{prop:two-secant-plane}
Suppose $p$ is odd and $S$ is not a cylinder. For each $P\in S$
there is a rational projective plane $M_P$ not containing $P$ such
that every other endpoint of a two-secant through $P$ belongs to
$M_P$.
\end{proposition}
\begin{proof}
Take rational coordinates with $P=e_0$. Since $P$ is fundamental,
the reduced complementary form has an expansion
\[
 G=z_0^{p-2}H+z_0^{p-3}K+
        \text{terms of lower $z_0$-degree},
\]
where $H$ is the local heavy-plane product of degree $p-1$ and $K$
has degree $p$. For each heavy plane through $P$, let its equation
be $\ell_i(z_1,z_2,z_3)$ and its weight be $r_i$. Lemma
\ref{lem:heavy-truncated-power}, on taking the coefficient of
$z_0^{p-3}$, gives
\[
 K\in(\ell_i^{r_i})+(z_1^p,z_2^p,z_3^p).
\]
The heavy normals span three dimensions: otherwise all heavy planes
would contain a common affine line through $P$, which the radial
identity would make a full line of $S$. By Lemma
\ref{lem:degree-p-gluing}, $K=B^p+LH$ for rational linear forms.

The fundamental equations imply that all ordinary partial
derivatives of $G$ of total order $p-2$ vanish at every $Q\in S$.
Indeed, for $|\beta|=p-2$, the coefficient of $z^\beta$ in
$D_Q^{p-1}G$ is $(p-1)!\,\partial^\beta G(Q)/\beta!$.
Therefore the gradient of
\[
 \frac{1}{(p-3)!}\partial_0^{p-3}G
       =K+(p-2)z_0H=B^p+(L-2z_0)H
\]
vanishes at all these points.

If $PQ$ is a two-secant, the radial identity says that its direction
belongs to exactly one heavy plane through $P$, and that this plane
has weight one. Thus $H(Q)=0$ and $\nabla H(Q)\ne0$. Evaluating the
last gradient at $Q$ gives $L(Q)-2Q_0=0$. Hence one can take
\[
 M_P=\{L(z_1,z_2,z_3)-2z_0=0\}.
\]
It is rational and misses $P$, since $p$ is odd.
\end{proof}

\section{The harmonic identity on every secant}\label{sec:harmonic}

We first make the normalization of the endpoint equation intrinsic.
Fix the original affine functional $w$. For $P\in S$, use coordinates
$z=sP+y$ with $y\in\ker w$ and $s=w(z)$, and take the truncated
representative of $G$. Its top two coefficients are
\[
 G_P(s,y)=s^{p-2}H_P(y)+s^{p-3}K_P(y)+\cdots,
 \qquad K_P=B_P^p+\ell_P H_P.
\]
Define the linear form on $V$ by
\begin{equation}\label{eq:canonical-lambda}
 \lambda_P(sP+y)=-2s+\ell_P(y).
\end{equation}
This is independent of the chosen basis of $\ker w$. Indeed, $H_P$
has degree $p-1$ and transforms as an ordinary form; $K_P$ has degree
$p$ and changes under truncation only by a pure $p$th power. The
linear quotient $\ell_P$ is unique: if $(\ell_P-\ell'_P)H_P$ were a nonzero
$p$th power, its linear-factor multiplicities would all be divisible
by $p$. But $H_P$ has at least three distinct heavy factors, with
positive multiplicities below $p$, and multiplication by one linear
form cannot make all of them divisible by $p$. Thus $\ell_P=\ell'_P$.
This proves that \eqref{eq:canonical-lambda} is one fixed covector
for each $P$, with $\lambda_P(P)=-2$, and its kernel is $M_P$.
In particular, the line-adapted calculations below use the same
forms that enter the final endomorphism.

\begin{proposition}\label{prop:secant-harmonic-formula}
Let an affine line meet a hypothetical noncylinder in
$P_i=P_0+t_i v$, $1\le i\le m$, where $2\le m<p$ and $v$ is a
fixed nonzero vector at infinity. Then
\[
 \lambda_{P_i}(v)=\frac{2}{m-1}
                    \sum_{j\ne i}\frac{1}{t_j-t_i}.
\]
In particular, if $P,Q$ lie on a three-secant, then
\[
 \lambda_P(Q)\lambda_Q(P)=1.
\]
If they lie on a two-secant, both factors are zero.
\end{proposition}
\begin{proof}
Choose coordinates $(s,t,y,z)$ with the indicated line given by
$y=z=0$, taken in the canonical family of
Section~\ref{sec:harmonic}: the first basis vector is a point of the
line with $w=1$, and the remaining basis vectors, including $v$, lie
in $\ker w$. The translation of the parameter below stays in this
family, so the linear form $L$ computed there is the canonical
$\ell_{P_i}$ of \eqref{eq:canonical-lambda}. We use the following binary calculation. A reduced binary
form $F(s,t)$ of degree $2N-m$ determines a degree-$m$ form
\[
 C_F(a,b)=\sum_{u+v=m}(-1)^u\binom{m}{u}
               [s^{N-u}t^{N-v}]F\,a^v b^u.
\]
The output $D_{(a,b)}^NF$ has degree $N-m<p$ and is annihilated
by $D_{(a,b)}$, since every $p$th ordinary directional derivative
is zero in characteristic $p$. It is therefore a scalar multiple
of $(bs-at)^{N-m}$. Comparing the coefficient of $t^{N-m}$ gives
\[
 D_{(a,b)}^N F
       =(-1)^mN!\,C_F(a,b)(bs-at)^{N-m}.
\]
All binomial coefficients used here are nonzero because $m<p$.
Since a reduced binary form of degree $2N-m$ has only the monomials
$s^{N-u}t^{N-v}$ with $u+v=m$, each of which appears in $C_F$ with a
nonzero coefficient, the map $F\mapsto C_F$ is injective on reduced
binary forms of degree $2N-m$.
The coefficient comparison can first be made for $a\ne0$ and then
extended as a polynomial identity in $a,b$.

Expand $G$ in the transverse variables $y,z$. Its coefficient of
transverse degree $k$ is binary of degree $2N-1-k$. Since all $m$
points on the line are fundamental, the associated binary form of
degree $k+1$ vanishes at them. For $k<m-1$, this associated form has
more zeros than its degree and is zero, so by injectivity the
binary coefficients themselves are zero. For $k=m-1$, the associated
form of each binary coefficient is a scalar multiple of
\[
 C(s,t)=\prod_{j=1}^m(t-t_j s).
\]
Thus the first possible transverse term is
$F(s,t)J(y,z)$, where $C_F$ is proportional to $C$.
It is nonzero: the local heavy-plane product at any $P_i$ has
order exactly $m-1$ in the direction of this line, by the radial
identity.

Translate the parameter so that $t_i=0$ and $P_i=e_0$. Write
$c_1,c_2$ for the coefficients of $s^{m-1}t,s^{m-2}t^2$ in $C$.
Then
\[
 \frac{c_2}{c_1}=-\sum_{j\ne i}\frac1{t_j}.
\]
If $f_1,f_2$ are the coefficients of
$s^{N-1}t^{N+1-m}$ and $s^{N-2}t^{N+2-m}$ in $F$, the binary
formula gives, for the same nonzero scalar $\kappa$,
\[
 c_1=-\kappa m f_1,\qquad c_2=\kappa\binom m2 f_2.
\]
Consequently,
\[
 \frac{f_2}{f_1}=\frac{2}{m-1}
                          \sum_{j\ne i}\frac1{t_j}.
\]
In the relation $K=B^p+LH$ from the endpoint-plane proof, compare
the terms of transverse degree $m-1$. The $p$th-power term has no
such term, since $1\le m-1<p$. Also
$H$ has no transverse terms below $m-1$, so the transverse part of
$L$ contributes only at larger degree. At degree $m-1$, only
$L(v)t$ times the first transverse block of $H$ remains, giving
$L(v)=f_2/f_1$. Since $\lambda_{P_i}=L-2s$, this proves the formula.

For three points with parameters $0,a,b$, it gives
\[
 \lambda_{P_0}(P_a)=\frac{a-b}{b},\qquad
 \lambda_{P_a}(P_0)=\frac{b}{a-b},
\]
whose product is one. The two-secant statement follows either from
the formula or from Proposition~\ref{prop:two-secant-plane}.
\end{proof}

\section{A light plane section closes the argument}\label{sec:completion}
\begin{proof}[Proof of Theorem~\ref{thm:prime-strong-cylinder}]
For $p=2,3$, the total positive heavy weight through any point
of $S$ is $p-1\le2$.  Thus there are at most two heavy planes
through that point.  They have a common affine line through
the point, and its radial heavy weight is $p-1$.  This line
is full, so Lemma~\ref{lem:full-line} applies.

Suppose $p\ge5$, and assume for contradiction that $S$ is not a
cylinder.  Lemmas~\ref{lem:full-line} and
\ref{lem:almost-full-line} exclude lines containing $p$ or
$p-1$ points of $S$.

Choose a plane $H$ with exactly $p$ points of $S$.  Such a
plane exists: through any point of $S$, at most $p-1$ planes
are heavy, whereas there are $p^2+p+1$ planes in total, and
every remaining plane through that point has intersection
size $p$.  Put $T=S\cap H$, and let $\widetilde H$ be the
three-dimensional vector subspace corresponding to the
projective closure of $H$.

For each $Q\in T$, use the canonical linear form $\lambda_Q$
from Proposition~\ref{prop:secant-harmonic-formula}, with
$\lambda_Q(Q)=-2$ for the original affine representative
$Q=(1,q)$.  Define the linear operator
\[
 A_H=\sum_{Q\in T}Q\,\lambda_Q|_{\widetilde H}
       \in\operatorname{End}(\widetilde H).
\]
Every line in $H$ has the same intersection with $T$ as with
$S$.  The harmonic formula therefore applies to every secant
of $T$.  Its line-sum identity is
\[
 \sum_{Q\in T\cap\ell}Q\lambda_Q(P)=-mP
 \quad\text{if }P\in T\cap\ell,\quad |T\cap\ell|=m.
\]
For completeness, writing the line points as $u+t_i v$ and
$a_i=\lambda_{u+t_i v}(v)$, the harmonic formula gives
\[
 \sum_i a_i=0,\qquad
 \sum_i t_i a_i=-m,\qquad
 \sum_i t_i^2a_i=-2\sum_i t_i;
\]
these three equalities give the displayed vector identity.
All denominators are invertible because $2\le m<p$.

Let $d_T(P)$ count the secants of $T$ through $P$.  Partition
$T\setminus\{P\}$ by these lines and remove the summand
$P\lambda_P(P)=-2P$ on each line.  Since $|T|=p$, this gives
\[
 A_H(P)
 =-2P+\sum_{\ell\ni P\text{ secant}}(2-|T\cap\ell|)P
 =(d_T(P)-p-1)P=(d_T(P)-1)P.
\]
Thus all $p$ points of $T$ are rational eigenvectors of the
same three-dimensional linear operator.

If $A_H$ is nonscalar, its rational projective eigenvectors
are contained in a line together with at most one further
point, or in at most three points.  Indeed, eigenspaces for
distinct eigenvalues have total vector-space dimension at
most three; nontrivial Jordan blocks add no eigenvectors.
Since $p\ge5$, at least $p-1$ points of $T$ must then be
collinear, contrary to Lemma~\ref{lem:almost-full-line}.

Therefore $A_H=\mu I$.  Its trace is
\[
 3\mu=\operatorname{tr}A_H
      =\sum_{Q\in T}\lambda_Q(Q)=-2p=0.
\]
As $p\ge5$, this implies $A_H=0$.  Consequently
$d_T(P)\equiv1\pmod p$.  But
$1\le d_T(P)\le p-1$, so $d_T(P)=1$.  All $p$ points of $T$
lie on its unique secant through $P$, giving a full affine
line and the final contradiction.
\end{proof}

\paragraph{Scope of the conclusion.}
The contradiction above covers every prime $p$. The small primes
$2,3$ use only the elementary local weights; every prime $p\ge5$
uses the same harmonic identity and light-plane trace argument.
There are no residual parameter branches in this proposed proof.
Its mathematical status remains the complete proof candidate stated
on the first page, pending external review.

\paragraph{Declaration of use of generative AI and AI-assisted technologies.}
This research was assisted by AI tools, in parts of mathematical exploration and drafting. The tools used: {ChatGPT Astra and Claude Opus 5.5}.
The author reviewed and edited all content and takes full responsibility for
the manuscript.

\end{document}